\documentstyle[12pt]{article}
\title{Nonclassical Lagrangian Dynamics and Potential Maps}
\author{Constantin Udri\c ste \\ {} \\
University Politehnica of Bucharest \\
Department of Mathematics I \\
Splaiul Independen\c tei 313 \\
77206 Bucharest, Romania\\
email:udriste@mathem.pub.ro}

\def\interior{{\setlength{\unitlength}{1em}\begin{picture}(0.85,1)
	\put(.1,0){\line(1,0){.55}}\put(.65,0){\line(0,1){.55}}\end{picture}}}

\date{}
\begin{document}
\maketitle

\newcommand{\ty}{\infty}
\newcommand{\ov}{\over}
\newcommand{\di}{\displaystyle}
\newcommand{\si}{\sigma}
\newcommand{\na}{\nabla}
\newcommand{\pa}{\partial}
\newcommand{\la}{\lambda}
\newcommand{\al}{\alpha}
\newcommand{\ep}{\varepsilon}
\newcommand{\va}{\varphi}
\newcommand{\ld}{\ldots}
\newcommand{\noa}{\noalign{\medskip}}

\begin{abstract}
Section 1 refines the theory of harmonic and potential maps. Section 2 defines
a generalized Lorentz world-force law and shows that any PDEs system of order
one generates such a law in suitable geometrical structure. In other words,
the solutions of any PDEs system of order one are harmonic or potential maps,
if we use semi-Riemann-Lagrange structures. Section 3 formulates open
problems regarding the geometry of semi-Riemann manifolds $(J^1(T,M), S_1)$,
$(J^2(T,M), S_2)$, and shows that the Lorentz-Udriste world-force
law is equivalent to covariant Hamilton PDEs on $(J^1(T,M), S_1)$.
\end{abstract}



\section{Harmonic Maps and Potential Maps}

All maps throughout the paper are smooth, while manifolds are real,
finite-dimensional, Hausdorff, second-countable and connected.

Let $(T,h)$ and $(M,g)$ be semi-Riemann manifolds of dimensions $p$ and $n$.
Hereafter we shall assume that the manifold $T$ is oriented. Greek (Latin)
letters will be used for indexing the components of geometrical objects
attached to the manifold $T$ (manifold $M$). Local coordinates will be
written
$$
t = (t^\al), \quad \al = 1, \ld, p
$$
$$
x = (x^i), \quad i = 1, \ld, n,
$$
and the components of the corresponding metric tensor and Christoffel symbols
will be denoted by $h_{\al \beta}, g_{ij}$, $H^\al_{\beta \gamma}$, $G^i_{jk}$.
Indices of tensors or distinguished tensors will be rised and lowered in the
usual fashion.

Let $\va:T \to M$, $\va (t) = x$, $x^i = x^i (t^\al)$ be a $C^\ty$ map
(parametrized sheet). We set
$$
x^i_\al = {\pa x^i \ov \pa t^\al}, \; x^i_{\al\beta} = {\pa^2 x^i \ov
\pa t^\al \pa t^\beta} - H^\gamma_{\al\beta} x^i_\gamma + G^i_{jk} x^j_\al
x^k_\beta. \leqno (1)
$$
Then $x^i_{\al}, \; x^i_{\al\beta}$ transform like tensors under coordinate
transformations $t \to \bar t$, $x \to \bar x$. In the sequel
$x^i_{\al}, \; x^i_{\al\beta}$ will be interpreted like distinguished tensors.

The canonical form of the energy density $E(\va)$ of the map $\va$ is defined
by
$$
E_0 (\va) (t) = {1 \ov 2} h^{\al\beta} (t) g_{ij} (x(t)) x^i_\al (t)
x^j_\beta (t).
$$
For a relatively compact domain $\Omega \subset T$, we define the energy
$$
E_0 (\va, \Omega) = \int_\Omega E_0 (\va) (t) dv_h,
$$
where $dv_h = \sqrt{|h|} dt^1 \wedge \ld \wedge dt^p$ denotes the volume
element induced by the semi-Riemann metric $h$. A map $\va$ is called
{\it harmonic map} if it is a critical point of the energy functional $E_0$,
i.e., an extremal of the Lagrangian
$$
L = E_0 (\va) (t) \sqrt{|h|},
$$
for all compactly supported variations (It should be remarked that every
$C^2$ harmonic map is automatically $C^\ty$). The {\it harmonic map equation}
is a system of nonlinear PDEs of second order
(generalized Laplace equations) and is expressed in local coordinates as
$$
\tau (\va)^i = h^{\al\beta} x^i_{\alpha\beta} = 0. \leqno (2)
$$
The quantity $\tau (\va)^i$ defines a section of the pull-back bundle
$\va^{-1} TM$ of the tangent bundle $TM$ of the manifold $M$ along $\va$, and
is called the {\it tension field} of $\va$.

The product manifold $T \times M$ is coordinated by $(t^\al, x^i)$. The
1-order jet manifold $J^1 (T,M)$, i.e., the configuration bundle, is endowed
with the adapted coordinates $(t^\al, x^i, x^i_\al)$. The distinguished tensors fields
and other distinguished geometrical objects
on $T \times M$ are introduced [4] using the jet bundle $J^1 (T,M)$.

Let $X^i_\al (t,x)$ be a given $C^\ty$ distinguished tensor field on $T \times M$ and
$c(t,x)$ be a given $C^\ty$ real function on $T\times M$. The general energy
density $E(\va)$ of the map $\va$ is defined by
$$
E(\va (t)) = {1 \ov 2} h^{\al\beta} (t) g_{ij} (x(t)) x^i_\al (t) x^j_\beta (t)
- h^{\al\beta} (t) g_{ij} (x(t)) x^i_\al (t) X^j_\beta (t, x(t)) + c(t,x).
$$
Of course $E(\va)$ is a perfect square iff
$$
c = {1 \ov 2} h^{\al\beta}(t) g_{ij} (x(t)) X^i_\al (t, x(t))
X^j_\beta (t, x(t)).
$$
Similarly, for a relatively compact domain $\Omega \subset T$, we define the
energy
$$
E(\va; \Omega) = \int_\Omega E(\va) (t) dv_h.
$$
A map $\va$ is called {\it potential map} if it is a critical point of the
energy functional $E$, i.e., an extremal of the Lagrangian
$$
L = E(\va)(t) \sqrt{|h|},
$$
for all compactly supported variations. The {\it potential map equation} is a system
of nonlinear PDEs (generalized Poisson equations)
and is expressed locally by
$$
\tau (\va)^i = h^{\al\beta} x^i_{\al\beta} = g^{ij} {\pa c \ov \pa x^j} +
h^{\al\beta} (\na_k X^i_\beta - g_{kj} g^{il} \na_l X^j_\beta)x^k_\al +
h^{\al\beta} D_\al X^i_\beta, \leqno (3)
$$
where $D$ is the Levy-Civita connection of $(T,h)$ and $\na$ is the
Levy-Civita connection of $(M,g)$. Explicitly, we have
$$
\na_j X^i_\al = {\pa X^i_\al \ov \pa x^j} + G^i_{jk} X^k_\al, \quad
D_\beta X^i_\al = {\pa X^i_\al \ov \pa t^\beta} - H^\gamma_{\beta\al} X^i_\gamma
\leqno (4)
$$
$$
F_j{}^i{}_\al = \na_j X^i_\al - g_{hj} g^{ik} \na_k X^h_\al, \leqno (5)
$$
$$
{\pa g_{ij} \ov \pa x^k} = G^h_{ki} g_{hj} + G^h_{kj} g_{hi}, \quad
{\pa h^{\al\beta} \ov \pa t^\gamma} = - H^\al_{\gamma\lambda} h^{\la \beta}
- H^\beta_{\gamma\la} h^{\al\la}. \leqno (6)
$$

\section{Lorentz-Udri\c ste World-Force Law}

In nonquantum relativity there are three basic laws for particles: the
Lorentz World-Force Law and two conservation laws [5]. Now we shall generalize
the Lorentz World-Force Law (see also [6], [11]).

{\bf Definition}. Let $F_\al =(F_j{}^i{}_\al)$ and $U_{\al\beta} =
(U^i_{\al\beta})$ be $C^\ty$ distinguished tensors on $T\times M$,
where $\omega_{ji\al} = g_{hi}F_j{}^h{}_\al$ is skew-symmetric
with respect to $j$ and $i$. Let $c(t,x)$
be a $C^\ty$ real function on $T \times M$. A $C^\ty$ map $\va: T \to M$ obeys the
{\it Lorentz-Udri\c ste World-Force Law} with respect to $F_\al$, $U_{\al\beta}$,
$c$ iff
$$
\tau (\va)^i = g^{ij} {\pa c \ov \pa x^j} + h^{\al\beta} F_j{}^i{}_\al
x^j_\beta + h^{\al\beta} U^i_{\al\beta}. \leqno (7)
$$

Now we show that the solutions of a system of PDEs of order one are potential
maps in a suitable geometrical structure. First we remark that a $C^\ty$
distinguished tensor field $X^i_\al (t,x)$ on $T \times M$ defines a family of
$p$-dimensional sheets as solutions of the PDEs system of order one
$$
x^i_\al = X^i_\al (t, x(t)), \leqno (8)
$$
if the complete integrability conditions
$$
{\pa X^i_\al \ov \pa t^\beta} + {\pa X^i_\al \ov \pa x^j} X^j_\beta =
{\pa X^i_\beta \ov \pa t^\al} + {\pa X^i_\beta \ov \pa x^j} X^j_\al
$$
are satisfied.

The distinguished tensor field $X^i_\al$ and semi-Riemann metrics $h$ and $g$
determine the {\it potential energy}
$$
f: T \times M \to R, \quad f = {1 \ov 2} h^{\al\beta} g_{ij} X^i_\al
X^j_\beta.
$$
The distinguished tensor field (family of $p$-dimensional sheets) $X^i_\al$ on
$(T \times M, h + g)$ is called:

1) {\it timelike}, if $f < 0$;

2) {\it nonspacelike or causal}, if $f \le 0$;

3) {\it null or lightlike}, if $f=0$;

4) {\it spacelike}, if $f>0$.

Let $X^i_\al$ be a distinguished tensor field of everywhere constant energy.
If $X^i_\al$ (the system (8)) has no critical point on $M$, then upon rescaling,
it may be supposed that
$f \in \left\{ -\di{1 \ov 2}, 0, \di{1 \ov 2} \right\}$.
Generally, $\cal E \subset M$ is the set of critical points of the distinguished
tensor  field $X^i_\al$, and this rescaling is possible only on $T \times (M \setminus \cal {E})$.

Using the operator (derivative along a solution of (8)),
$$
{\delta \ov \pa t^\beta}x^i_\al = x^i_{\al\beta} =
{\pa^2 x^i \ov \pa t^\al \pa t^\beta} -
H^\gamma_{\al\beta} x^i_\gamma + G^i_{jk} x^j_\al x^k_\beta
$$
we obtain the prolongation (system of PDEs of order two)
$$
x^i_{\al\beta} = D_\beta X^i_\al + (\na_j X^i_\al) x^j_\beta. \leqno (9)
$$

The distinguished tensor field $X^i_\al$, the metric $g$, and the connection $\na$
determine the external distinguished tensor field
$$
F_j{}^i{}_\al = \na_j X^i_\al - g^{ih} g_{kj} \na_h X^k_\al,
$$
which characterizes the {\it helicity} of the distinguished tensor field
$X^i_\al$.

First we write the PDEs system (9) in the equivalent form
$$
x^i_{\al\beta} = g^{ih} g_{kj} (\na_h X^k_\al) x^j_\beta + F_j{}^i{}_\alpha
x^j_\beta + D_\beta X^i_\al.
$$
Now we modify this PDEs system into
$$
x^i_{\al\beta} = g^{ih} g_{kj} (\na_h X^k_\al)X^j_\beta + F^i_{j\al} x^j_\beta
+ D_\beta X^i_\al. \leqno (10)
$$
Of course, the PDEs system (10) is still a prolongation of the PDEs system
(8).

Taking the trace of (10) with respect to $h^{\al\beta}$ we obtain that any
solution of PDEs system (8) is also a solution of the PDEs system
$$
h^{\al\beta} x^i_{\al\beta} = g^{ih} h^{\al\beta} g_{kj} (\na_h X^k_\al)
X^j_\beta + h^{\al\beta} F_j{}^i{}_\al x^j_\beta + h^{\al\beta} D_\beta X^i_\al.
\leqno (11)
$$
(generalized Poisson equations).

{\bf Theorem}. {\it The PDEs system (11) is a prolongation of the PDEs
system (8)}.

If $F_j{}^i{}_\al = 0$, then the PDEs system (11) reduces to
$$
h^{\al\beta} x^i_{\al\beta} = g^{ih} h^{\al\beta} g_{kj} (\na_h X^k_\al)X^j_\beta
+ h^{\al \beta} D_\beta X^i_\al. \leqno (12)
$$

The first term in the second hand member of the PDEs systems (11) or (12)
is $(grad \: f)^i$. Consequently, choosing the metrics $h$ and $g$ such
that $f \in \left\{ -\di{1 \ov 2}, 0, \di{1 \ov 2} \right\}$, then the
preceding PDEs systems reduce to
$$
h^{\al\beta} x^i_{\al\beta} = h^{\al\beta} F_j{}^i{}_\al x^j_\beta +
h^{\al\beta} D_\beta X^i_\al \leqno (11')
$$
$$
h^{\al\beta} x^i_{\al\beta} = h^{\al\beta} D_\beta X^i_\al. \leqno (12')
$$

{\bf Theorem}. 1) {\it The solutions of PDEs system (11) are the extremals
of the Lagrangian}
$$
\begin{array}{lcl}
L &=& \di{1 \ov 2} h^{\al\beta} g_{ij} (x^i_\al - X^i_\al)(x^j_\beta - X^j_\beta)
\sqrt{|h|} = \\ \noa
&=& \left( \di{1 \ov 2} h^{\al\beta} g_{ij} x^i_\al x^j_\beta -
h^{\al\beta} g_{ij} x^i_\al X^j_\beta + f\right) \sqrt{|h|}. \end{array}
$$

2) {\it The solution PDEs system (12) are the extremals of the Lagrangian}
$$
L = \left( {1 \ov 2} h^{\al\beta} g_{ij} x^i_\al x^j_\beta + f \right)
\sqrt{|h|}.
$$

3) {\it If the Lagrangians $L$ are independent of the variable $t$, then the
PDEs systems (11) or (12) are conservative, the energy-impulse tensor field
being}
$$
T^\al{}_\beta = x^i_\beta {\pa L \ov \pa x^i_\al} - L \delta^\al_\beta.
$$

4) {\it Both Lagrangians produce the same Hamiltonian}
$$
H = \left( {1 \ov 2} h^{\al\beta} g_{ij} x^i_\al x^j_\beta - f \right)
\sqrt{|h|}.
$$

{\bf Proof}. 1) and 2) If we write $L = E \sqrt{|h|}$, where $E$ is the energy
density, then the Euler-Lagrange equations of extremals
$$
{\pa L \ov \pa x^k} - {\pa \ov \pa t^\al} {\pa L \ov \pa x^k_\al} = 0
$$
can be written
$$
{\pa E \ov \pa x^k} - {\pa \ov \pa t^\al}{\pa E \ov \pa x^k_\al} -
H^\gamma_{\gamma \al} {\pa E \ov \pa x^k_\al} = 0. \leqno (13)
$$
We compute
$$
\begin{array}{lcl}
\di{\pa E \ov \pa x^k} &=& \di{1 \ov 2} h^{\al\beta} \di{\pa g_{ij} \ov
\pa x^k} x^i_\al x^j_\beta - h^{\al\beta} \di{\pa g_{ij} \ov \pa x^k}
x^i_\al X^j_\beta + \\ \noa
&+& \di{1 \ov 2} h^{\al\beta} \di{\pa g_{ij} \ov \pa x^k} X^i_\al X^j_\beta -
h^{\al\beta} g_{ij} x^i_\al \di{\pa X^j_\beta \ov \pa x^k} +
h^{\al\beta} g_{ij} \di{\pa X^i_\al \ov \pa x^k} X^j_\beta, \end{array}
$$

\hspace{1 cm} $\di{\pa E \ov \pa x^k_\al} = h^{\al\beta}
g_{kj} x^j_\beta - h^{\al\beta} g_{kj} X^j_\beta,$

$$
\begin{array}{lcl}
- \di{\pa \ov \pa t^\al}\di{\pa E \ov \pa x^k_\al} &=& - \di{\pa h^{\al\beta}
\ov \pa t^\al} g_{kj} x^j_\beta - h^{\al\beta} \di{\pa g_{kj} \ov \pa x^l}
x^l_\al x^j_\beta - h^{\al\beta} g_{kj} \di{\pa^2x^j \ov\pa t^\al \pa t^\beta} + \\ \noa
&+& \di{\pa h^{\al\beta} \ov \pa t^\al} g_{kj} X^j_\beta +
h^{\al\beta} \di{\pa g_{kj} \ov \pa x^l} x^l_\al X^j_\beta +
h^{\al\beta} g_{kj} \left(\di{\pa X^j_\beta \ov \pa t^\al} +
\di{\pa X^j_\beta \ov \pa x^l} x^l_\al \right). \end{array}
$$
We replace in (13) taking into account the formulas (1), (4) and (6). We
find
$$
\begin{array}{lcl}
h^{\al\beta} g_{kj} x^j_{\al\beta} &=& h^{\al\beta} g_{ij} (\na_k X^i_\al)
X^j_\beta + h^{\al\beta} g_{kj} (\na_l X^j_\beta)x^l_\al - \\ \noa
&-& h^{\al\beta} g_{ij} x^i_\al \na_k X^j_\beta +  h^{\al\beta} g_{kj}
D_\al X^j_\beta. \end{array}
$$

Transvecting by $g^{hk}$ and using the formula (5), we obtain
$$
h^{\al\beta} x^i_{\al\beta} = g^{ik}  h^{\al\beta} g_{lj} (\na_k X^l_\al)
X^j_\beta +  h^{\al\beta} F_j{}^i{}_\al x^j_\beta +  h^{\al\beta} D_\al X^i_\beta.
$$

3) Taking into account the Euler-Lagrange equations, we have
$$
{\pa T^\al_\beta \ov \pa t^\al} = {\pa^2 x^i \ov \pa t^\al \pa t^\beta}
{\pa L \ov \pa x^i_\al} + x^i_\beta {\pa^2 L \ov \pa t^\al \pa x^i_\al}
+ x^i_\beta {\pa^2 L \ov \pa x^j \pa x^i_\al} x^j_\al +
$$
$$
+ x^i_\beta {\pa^2 L \ov \pa x^i_\al \pa x^j_\gamma}
{\pa^2 x^j \ov \pa t^\gamma \pa t^\al}
- {\pa L \ov \pa t^\al} \delta^\al_\beta - {\pa L \ov \pa x^j} x^j_\beta -
{\pa L \ov \pa x^j_\gamma} {\pa^2 x^j \ov \pa t^\gamma \pa t^\al}
\delta^\al_\beta = - {\pa L \ov \pa t^\beta}.
$$

{\bf Open problem}. Determine the general expression of the energy-impulse
tensor field as object on $J^1(T,M)$, and compute its divergence.

4) We use the formula
$$
H = x^i_\al {\pa L \ov \pa x^i_\al} - L.
$$

{\bf Corollary}. {\it Every PDE generates a Lagrangian of order one via the
associated first order PDEs system and suitable metrics on the manifold of
independent variables and on the manifold of functions. In this sense the
solutions of the initial PDE are potential maps.}

{\bf Theorem (Lorentz-Udri\c ste World-Force Law)}.

1) {\it Every solution of the PDEs system (12) is a potential map on the
semi-Riemann manifold $(T \times M, h+g)$.}

2) {\it Every solution of the PDEs system (11) is a horizontal potential map of
the semi-Riemann-Lagrange manifold}
$$
\left( T \times M,\; h+g, \; N(^i_\al)_j = G^i_{jk} x^k_\al - F_j{}^i{}_\al, \;
M(^i_\al)_\beta = - H^\gamma_{\al\beta} x^i_\gamma \right).
$$

\section{Covariant Hamilton Field Theory \protect \\
(Covariant Hamilton equations)}

Let us show that the PDEs systems (11) and (12) induce on $J^1 (T,M)$
Hamilton PDEs systems.

Let $(T,h)$ be a semi-Riemann manifold with $p$ dimensions, and $(M,g)$ be a
semi-Riemann manifold with $n$ dimensions. Then
$(J^1(T,M), h+g + h^{-1} * g)$ is a semi-Riemann manifold with
$p+n+pn$ dimensions.

We denote by $X^i_\al$ a $C^\ty$ distinguished tensor field on $T \times M$, and by
$\omega_{ij\al}$ the distinguished 2-form associated to the
distinguished tensor field
$$
F_j{}^i{}_\al = \na_j X^i_\al - g^{ih} g_{kj} \na_h X^k_\al
$$
via the metric $g$, i.e., $\omega = \di{1 \ov 2} g \circ F$. Of course
$X^i_\al$, $F_j{}^i{}_\al$ are distinguished objects on $J^1 (T,M)$ globally defined.

If $(t^\al, x^i, x^i_\al)$ are the coordinates of a point in $J^1(T,M)$,
and $H^\al_{\beta\gamma}$, $G^i_{jk}$ are the components of the connection
induced by $h$ and $g$, respectively, then
$$
\left( {\delta \ov \delta t^\al} = {\pa \ov \pa t^\al} + H^\gamma_{\al\beta}
x^i_\gamma {\pa \ov \pa x^i_\beta}, \;
{\delta \ov \delta x^i} = {\pa \ov \pa x^i} - G^h_{ik} x^k_\al
{\pa \ov \pa x^h_\al}, \; {\pa \ov \pa x^i_\al}\right),
$$
$$
\left( dt^\beta, dx^j, \; \delta x^j_\beta = dx^j_\beta - H^\gamma_{\beta\la} x^j_\gamma
dt^\la + G^j_{hk} x^h_\beta dx^k \right)
$$
are dual frames on $J^1 (T,M)$, i.e.,
$$
dt^\beta\left({\delta \ov \delta t^\al}\right) = \delta^\beta_\al, \quad
dt^\beta \left({\delta \ov \delta x^i}\right) = 0, \quad
dt^\beta \left({\pa \ov \pa x^i_\al}\right) = 0
$$
$$
dx^j\left({\delta \ov \delta t^\al}\right) = 0, \quad
dx^j  \left({\delta \ov \delta x^i}\right) = \delta^j_i, \quad
dx^j \left({\pa \ov \pa x^i_\al}\right) = 0
$$
$$
\delta x^j_\beta\left({\delta \ov \delta t^\al}\right) = 0, \quad
\delta x^j_\beta \left({\delta \ov \delta x^i}\right) = 0, \quad
\delta x^j_\beta \left({\pa \ov \pa x^i_\al}\right) = \delta^j_i \delta^\al_\beta.
$$

Using these frames, the induced Sasaki-like metric on $J^1(T,M)$ is
$$
S_1 = h_{\al\beta} dt^\al \otimes dt^\beta + g_{ij} dx^i \otimes dx^j +
h^{\al\beta} g_{ij} \delta x^i_\al \otimes \delta x^j_\beta.
$$

{\bf Open problems}.

1) The geometry of the semi-Riemann manifold $(J^1 (T,M),S_1)$, which is similar
to the geometry of the tangent bundle endowed with Sasaki metric, is now in
working by our research group [4]. As was shown here this geometry permits the
interpretation of solutions of PDEs systems of order one (8) as potential
maps. In this sense the solutions of every PDE of any order are extremals of a
Lagrangian of order one.

2) Study the geometry of the dual space of $(J^1(T,M), S_1)$.

3) Find a Sasaki-like $S_2$ metric on the jet bundle of order two and develop
the geometry of the semi-Riemann manifold $(J^2(T,M),S_2)$. In this manifold,
the PDEs of Mathematical Physics (of order two) appear like hypersurfaces.
Most of them are in fact algebraic hypersurfaces.

4) Study the geometry of the dual space of $(J^2 (T,M),S_2)$.

Recall that on a symplectic manifold $(Q,\Omega)$ of even dimension $q$, the
Hamiltonian vector field $X_{f_1}$ of a function $f_1 \in \cal F(Q)$ is
defined by
$$
X_{f_1}  \interior \Omega = df_1,
$$
and the Poisson bracket of $f_1, f_2$ is defined by
$$
\{f_1, f_2\} = \Omega (X_{f_1}, X_{f_2}).
$$
The polysymplectic analogue of a function is a $q$-form called {\it momentum
observable}. The Hamiltonian vector field $X_{f_1}$ of such a momentum
observable $f_1$ is defined by
$$
X_{f_1}  \interior \Omega = df_1,
$$
where $\Omega$ is the canonical $(q+2)$-form  on the appropriate dual of
$J^1 (T,M)$. Since $\Omega$ is nondegenerate, this uniquely defines $X_{f_1}$.
The Poisson bracket of two such $n$-forms $f_1, f_2$ is the $n$-form defined
by
$$
\{f_1, f_2\} = X_{f_1}  \interior (X_{f_2} \interior \Omega ).
$$
Of course $\{f_1, f_2\}$ is, up to the addition of exact terms, another
momentum observable.

{\bf Theorem}. {\it The PDEs system
$$
h^{\al\beta} x^i_{\al\beta} = g^{ih} h^{\al\beta} g_{jk} X^j_\beta \na_h
X^k_\al
$$
transfers in $J^1(T,M)$ as a covariant Hamilton PDEs system with respect to
the Hamiltonian (momentum observable)
$$
H = \left( {1 \ov 2} h^{\al\beta} g_{ij} x^i_\al x^j_\beta - f\right) dv_h
$$
and the non-degenerate distinguished polysymplectic $(p+2)$-form}
$$
\Omega = \Omega_\al \otimes dt^\al, \quad \Omega_\al = g_{ij} dx^i \wedge
\delta x^j_\al \wedge dv_h.
$$

{\bf Proof}. Let
$$
\theta = \theta_\al \otimes dt^\al, \quad \theta_\al = g_{ij} x^i_\al dx^j
\wedge dv_h
$$
be the distinguished Liouville $(p+1)$-form on $J^1 (T,M)$. It follows
$$
\Omega_\al = -d\theta_\al.
$$
We denote by
$$
X_H = X^\beta_H {\delta \ov \delta t^\beta}, \quad
X^\beta_H = u^{\beta l} {\delta \ov \delta x^l} + {\delta u^{\beta l} \ov
\pa t^\al}{\pa \ov \pa x^l_\al}
$$
the distinguished Hamiltonian object of the observable $H$. Imposing
$$
X^\al_H \interior \Omega_\al = dH,
$$
where
$$
dH = (h^{\al\beta} g_{ij} x^j_\beta \delta x^i_\al - h^{\al\beta} g_{ij}
X^j_\beta \na_k X^i_\al dx^k) \wedge dv_h,
$$
we find
$$
g_{ij} u^{\al i} \delta x^j_\al - g_{ij} {\delta u^{\al j}\ov \pa t^\al}
dx^i = h^{\al\beta} g_{ij} x^j_\beta \delta x^i_\al -
h^{\al\beta} g_{ij} X^j_\beta \na_k X^i_\al dx^k \; \hbox{modulo}\; dv_h.
$$
Consequently, it appears the Hamilton PDEs system
$$
\left\{ \begin{array}{l}
u^{\al i} = h^{\al \beta} x^i_\beta \\ \noa
\di{\delta u^{\al i} \ov \pa t^\al} = g^{hi} h^{\al\beta} g_{jk}
X^j_\beta (\na_h X^k_\al) \end{array} \right.
$$
(up to the addition of terms which are cancelled by the exterior multiplication
with $dv_h$).

{\bf Theorem}. {\it The PDEs system
$$
h^{\al\beta} x^i_{\al\beta} = g^{ih} h^{\al\beta} g_{kj} (\na_h X^k_\al)
X^j_\beta + h^{\al\beta} F_j{}^i{}_\al x^j_\beta + h^{\alpha\beta}
D_\beta X^i_\alpha
$$
transfers in $J^1(T,M)$ as a covariant Hamilton PDEs system with respect to
the Hamiltonian (momentum observable)
$$
H = \left( {1 \ov 2} h^{\al\beta} g_{ij} x^i_\al x^j_\beta - f\right) dv_h
$$
and the non-degenerate distinguished polysymplectic $(p+2)$-form}
$$
\Omega = \Omega_\al \otimes dt^\al, \quad \Omega_\al =
(g_{ij} dx^i \wedge \delta x^j_\al + \omega_{ij\al} dx^i \wedge dx^j
+ g_{ij} (D_\beta X^i_\alpha) dt^\beta \wedge dx^j) \wedge dv_h.
$$

{\bf Proof}. Let
$$
\theta = \theta_\al \otimes dt^\al, \quad \theta_\al = (g_{ij} x^i_\al dx^j
- g_{ij} X^i_\al dx^j) \wedge dv_h
$$
be the distinguished Liouville $(p+1)$-form on $J^1(T,M)$. It follows
$$
\Omega_\al = -d\theta_\al
$$
(of course the term containing $dt^\beta$ disappears by exterior multiplication
with $dv_h$). We denote by
$$
X_H = X^\beta_H {\delta \ov \delta t^\beta}, \quad
X^\beta_H = h^{\beta\gamma} {\delta \ov \delta t^\gamma} +
u^{\beta l} {\delta \ov \delta x^l} + {\delta u^{\beta l} \ov
\pa t^\al}{\pa \ov \pa x^l_\al}
$$
the distinguished Hamiltonian object of the observable $H$. Imposing
$$
X^\al_H \interior \Omega_\al = dH,
$$
where
$$
dH = (h^{\al\beta} g_{ij} x^j_\beta \delta x^i_\al - h^{\al\beta} g_{ij}
X^j_\beta (\na_k X^i_\al) dx^k) \wedge dv_h,
$$
we find
$$
(g_{ij} u^{\al i} \delta x^j_\al - g_{ij} {\delta u^{\al j}\ov \pa t^\al}
dx^i + 2\omega_{ij\al} u^{\al i} dx^j + h^{\alpha\beta}g_{ij} (D_\beta X^i_\alpha)
dx^j) \wedge dv_h = dH.
$$
Consequently, it appears the Hamilton PDEs system
$$
\left\{ \begin{array}{l}
u^{\al i} = h^{\al \beta} x^i_\beta \\ \noa
\di{\delta u^{\al i} \ov \pa t^\al} = g^{hi} h^{\al\beta} g_{jk} X^j_\beta
(\na_h X^k_\al) + 2g^{hi} \omega_{j h \al} u^{\al j}
+ h^{\alpha\beta} D_\beta X^i_\alpha \end{array} \right.
$$
(up to the addition of terms which are cancelled by the exterior multiplication
with $dv_h$).

{\bf Example}. The $C^\ty$ vector fields $\xi_\al$ on the manifold $M$ and
the 1-forms $A^\al$ on the manifold $T$ satisfying
$$
[\xi_\al, \xi_\beta] = C^\gamma_{\al\beta} \xi_\gamma, \quad
C^\gamma{}_{\al\beta} = \hbox{constants},
$$
$$
{\pa A^\al_\beta \ov \pa t^\gamma} - {\pa A^\al_\gamma \ov \pa t^\beta}
= C^\al_{\la\delta} A^\la_\beta A^\delta_\gamma
$$
determine a continuous group of transformations via the PDEs
$$
x^i_\al = \xi^i_\beta (x(t)) A^\beta_\al (t).
$$

Conversely, if $x^i = x^i (t^\al, y^j)$ are solutions of a completely
integrable system of PDEs of the preceding form, where the $A's$ and $\xi's$
satisfy the conditions stated above, such that for values $t^\al_0$ of
$t's$ the determinant of the $A's$ is not zero and
$$
x^i (t^\al_0, y^j) = y^j,
$$
then $x^i = x^i (t^\al, y^j)$ define a continuous group of transformations.

Using a semi-Riemann metric $h$ on the manifold $T$, a semi-Riemann metric $g$
on the manifold $M$, then the maps determining a continuous group of
transformations appear like extremals (potential maps) of the Lagrangian
$$
L = {1 \ov 2} h^{\al\beta} g_{ij} (x^i_\al - \xi^i_\la A^\la_\al)
(x^j_\beta - \xi^j_\mu A^\mu_\beta) \sqrt{|h|}.
$$

\end{document}